\documentstyle{amsppt}
\magnification=1200
\hsize=150truemm
\vsize=224.4truemm
\hoffset=4.8truemm
\voffset=12truemm

\NoRunningHeads

\define\C{{\bold C}}
 
 
\let\thm\proclaim
\let\fthm\endproclaim

\define\bo{\partial  }

\define\De{\Delta}
\define\Dn{\Delta_n}

\newcount\tagno
\newcount\secno
\newcount\subsecno
\newcount\stno
\global\subsecno=1
\global\tagno=0
\define\ntag{\global\advance\tagno by 1\tag{\the\tagno}}

\define\sta{\ 
{\the\secno}.\the\stno
\global\advance\stno by 1}

\define\stas{\the\stno
\global\advance\stno by 1}

\define\sect{\global\advance\secno by 1
\global\subsecno=1\global\stno=1\
{\the\secno}. }

\def\nom#1{\edef#1{{\the\secno}.\the\stno}}
\def\inom#1{\edef#1{\the\stno}}
\def\eqnom#1{\edef#1{(\the\tagno)}}

\newcount\refno
\global\refno=0

\def\nextref#1{
      \global\advance\refno by 1
      \xdef#1{\the\refno}}

\def\bref {\ref\global\advance\refno by 1\key{\the\refno}}

\nextref\BRO
\nextref\BRU 
\nextref\DUV
\nextref\KLE
\nextref\MCQ
\nextref\BLO
\nextref\MAT
\nextref\PAU
\nextref\SIK

\topmatter

\title 
Sur le lemme de Brody
 \endtitle

\author  Julien Duval
\endauthor

\abstract\nofrills{\smc R\'esum\'e. }\ \ Le lemme de Brody est un outil de base
en hyperbolicit\'e complexe. On en pr\'esente une version pr\'ecisant la localisation d'une
courbe enti\`ere issue d'une suite divergente de disques holomorphes. Elle donne aussi
une  caract\'erisation de l'hyperbolicit\'e en termes d'in\'egalit\'e isop\'erim\'etrique.
 
 \endabstract

\endtopmatter 

\document

\subhead 0. Introduction \endsubhead
\stno=1

\null
\footnote""{Mots-cl\'es : hyperbolicit\'e complexe, courbe enti\`ere, in\'egalit\'e isop\'erim\'etrique}
\footnote""{Classification AMS : 32Q45, 32U40.}

Soit $X$ une vari\'et\'e compacte complexe, munie d'une m\'etrique hermitienne. Une suite divergente de
disques holomorphes dans $X$ produit par reparam\'etrage une courbe enti\`ere, i.e. une application holomorphe 
non constante de $\C$ dans $X$.
C'est le contenu du lemme de Brody :

\thm{Lemme [\BRO]} Soit $(f_n)$ une suite d'applications holomorphes du disque unit\'e $D$ de $\C$ dans $X$. On suppose
que $\Vert f'_n(0)\Vert$ tend vers l'infini. Alors il existe une suite $(r_n)$ de reparam\'etrages de $\C$
telle que $f_n \circ r_n$ converge localement uniform\'ement vers une courbe enti\`ere, apr\`es extraction.
\fthm

Rappelons que $X$ est {\it hyperbolique} (au sens de Kobayashi) si la pseudom\'etrique
 $$K(x,v)=\text{inf }\{\frac {1}{r} /\ \ \exists f:D\rightarrow X \text{ holomorphe avec } f(0)=x \text{ et } f'(0)=rv\} $$
est non d\'eg\'en\'er\'ee.
Le lemme de Brody permet de caract\'eriser les vari\'et\'es hyperboliques :

\thm{Corollaire} Une vari\'et\'e compacte complexe est hyperbolique si et seulement si elle ne contient pas de courbe enti\`ere.
\fthm

L'inconv\'enient principal de ce lemme est l'absence de contr\^ole sur la localisation de la courbe enti\`ere produite. Nous en
 pr\'esentons ici une version quantitative qui pallie ce d\'efaut. Pour cela on s'int\'eresse plut\^ot aux suites de disques
 holomorphes qui divergent du point de vue de l'aire. Il n'est pas difficile de supposer,
 quitte \`a r\'eduire un peu ces disques, que la longueur de leur bord devient n\'egligeable devant leur aire (voir [\BLO]).
On introduit donc tout naturellement la 

\null

\noindent
{\bf D\'efinition.} Un courant positif ferm\'e $T$ de masse 1 est dit {\it d'Ahlfors} s'il existe une suite $(\Dn)$ de disques
 holomorphes dans $X$
telle que long($\bo \Dn)=o(\text{aire}(\Dn))$ et
 $T=\lim \frac{[\Dn]}{\text{aire}(\Dn)}$.

\null

Ces courants apparaissent notamment dans la preuve par M. McQuillan ([\MCQ]) de la conjecture de Green-Griffiths pour
 certaines surfaces (voir aussi [\BRU],[\DUV],[\BLO],[\PAU]).

 Voici notre r\'esultat principal. Il localise la courbe
 enti\`ere produite \`a partir de la suite divergente de disques holomorphes l\`a o\`u l'aire s'accumule :

\thm{Th\'eor\`eme} Soit $T$ un courant d'Ahlfors dans $X$. On suppose que
$T$ charge un compact $K$. Alors il existe une courbe enti\`ere coupant $K$ sur un ensemble d'aire non nulle. \fthm

Une analyse de la d\'emonstration donne de plus que la courbe enti\`ere est contenue dans le support de $T$.

\null

Comme dans le lemme de Brody, cette courbe est obtenue par reparam\'etrage des disques $\Dn$ et passage \`a la limite,
cette fois au sens de Gromov. Pr\'ecis\'ement, si $\Dn=f_n(D)$, on construit une suite $(r_n)$ de reparam\'etrages de $\C$
 telle que :

i) $f_n \circ r_n$
soit d\'efinie asymptotiquement sur $\C$ (resp. $\C^*$),

ii) l'aire de l'image de $f_n \circ r_n$ soit uniform\'ement
major\'ee sur tout compact de $\C$ (resp. $\C^*$),

iii) l'image de $f_n\circ r_n$
 rencontre un voisinage de plus en plus petit
 de $K$ sur un ensemble d'aire au moins 1.

 Apr\`es extraction et du fait des bornes d'aire, $f_n\circ r_n$ converge localement uniform\'ement vers $f$
 hors d'un ensemble discret de $\C$ (resp. $\C^*$) et $f$ se prolonge holomorphiquement aux points de cet ensemble (voir
 par exemple
 [\SIK]).
Cette application (resp. son rel\`evement via l'exponentielle) est la courbe enti\`ere voulue.

\null

Voici comment on obtient ces reparam\'etrages. Notons pour simplifier $a_n$ l'aire de $\Dn$.
On produit d'abord, par un lemme de Besicovitch (voir [\MAT]),
 de petits disques disjoints dans $D$ en nombre de l'ordre de $a_n$ dont l'image par $f_n$ coupe
 un voisinage de $K$ sur un ensemble d'aire au moins 1. Cela vient du fait que $T$ charge $K$.
Comme ces disques sont disjoints et nombreux, on en d\'eduit une borne sur l'aire de leur image.  Ce sont les germes
 de la courbe enti\`ere. Puis on double ces disques toujours \`a la source. Un premier probl\`eme
 est de rester dans $D$. On l'assurera pour une bonne partie des disques car la longueur de $\bo \Dn$
devient n\'egligeable devant $a_n$. Un deuxi\`eme probl\`eme est de les conserver disjoints
en nombre suffisant. Si c'est possible on obtient une borne d'aire de l'image de ces disques doubl\'es. Sinon une
 analyse combinatoire de la r\'epartition des disques initiaux fournit des anneaux
disjoints de module arbitrairement grand et en nombre de l'ordre de $a_n$. Quand on peut r\'ep\'eter ind\'efiniment
 l'op\'eration de doublement on construit une image de $\C$. Sinon on produit une image de $\C^*$.

\null

Voici deux cons\'equences du th\'eor\`eme. La premi\`ere est une caract\'erisation de l'hyperbolicit\'e en termes 
d'in\'egalit\'e isop\'erim\'etrique, \`a la Gromov (voir [\KLE] pour  un \'enonc\'e analogue) :

\thm{Corollaire} Une vari\'et\'e compacte complexe est hyperbolique si et seulement si ses disques holomorphes v\'erifient une
in\'egalit\'e isop\'erim\'etrique lin\'eaire.
\fthm

Autrement dit, il existe une constante $C$ telle que tout disque holomorphe $\De$ v\'erifie 
 $$\text{aire}(\De)\leq C\text{ long}(\bo \De).$$ Cette propri\'et\'e est ind\'ependante de la m\'etrique hermitienne
choisie sur $X$. Quand elle n'est pas satisfaite on peut construire
un courant d'Ahlfors dans $X$, donc une courbe enti\`ere : $X$ n'est pas hyperbolique. Inversement, si $X$ n'est pas hyperbolique
on dispose d'une courbe enti\`ere. Le lemme d'Ahlfors (voir [\BRU]) nous donne alors une suite de disques concentriques
dans cette courbe dont la longueur de bord devient n\'egligeable devant l'aire : l'in\'egalit\'e isop\'erim\'etrique lin\'eaire
est viol\'ee.

\null

La seconde cons\'equence g\'en\'eralise un r\'esultat de [\DUV] :

\thm{Corollaire} Si un courant d'Ahlfors de $X$ charge un ensemble analytique, celui-ci contient
une courbe enti\`ere.
\fthm

En effet une courbe enti\`ere coupant un ensemble analytique sur un ensemble d'aire non nulle doit y \^etre contenue.

\null

Signalons que le th\'eor\`eme peut produire dans certains cas plusieurs courbes enti\`eres 
\`a partir d'une suite divergente de disques. Il suffit par exemple que le courant d'Ahlfors associ\'e charge 
un ensemble analytique sans \^etre enti\`erement port\'e par lui.

Les \'enonc\'es pr\'ec\'edents s'appliquent aussi aux courants d'Ahlfors limites d'une suite d'unions finies de disques holomorphes
 (voir [\DUV]).
Ceci permet de couvrir les courants issus de courbes enti\`eres apr\`es r\'egularisation \`a la Nevanlinna.

Enfin la souplesse de la m\'ethode employ\'ee rend (mutatis mutandis) ce th\'eor\`eme valide dans un cadre beaucoup
plus vaste, en fait d\`es qu'on dispose d'un th\'eor\`eme de compacit\'e \`a la Gromov :
 courbes pseudoholomorphes, applications harmoniques
minimisantes etc...

\null
Cet article r\'epond \`a des questions de M. McQuillan et M. Paun. Qu'ils en soient chaleureusement remerci\'es.

\newpage

\subhead 1. Pr\'eliminaires \endsubhead

\null

Voici l'\'enonc\'e de compacit\'e au coeur de ce qui va suivre.
Il remonte au moins \`a Gromov dans un contexte bien plus large (voir par exemple [\SIK]).
\thm {Th\'eor\`eme}
Soit $g_n:D \rightarrow X$ une suite d'applications holomorphes du disque unit\'e de $\C$ dans $X$.
 On suppose l'aire de $g_n(D)$ uniform\'ement born\'ee. Alors, quitte \`a extraire, 

 i) $g_n$ converge localement uniform\'ement vers $g$ sur $D$ priv\'e d'un nombre fini
de points d'explosion;

ii) $g$ se prolonge en une application holomorphe encore not\'e $g:D \rightarrow X$;

iii) en un point d'explosion $e$, il existe une suite de disques $d_n$ dans $D$ tendant vers $e$ telle que
$g_n(d_n)$ converge au sens de Hausdorff et en aire vers une union finie de courbes rationnelles (des bulles).

\fthm

Dans la suite on dira que $g_n$ tend vers $g$ {\it au sens de Gromov}. Les bulles sont en nombre fini du fait
de la borne d'aire. En effet une courbe rationnelle de $X$ a une aire minor\'ee a priori, par exemple par 1
quitte \`a normaliser la m\'etrique.

\null

En voici deux variantes \`a bord.

\thm{Variante 1} Soit $d_n$ une suite de disques convergeant vers 0 et $g_n:D-d_n \rightarrow X$
 une suite d'applications holomorphes, lisses jusqu'\`a $\bo d_n$.
On suppose que l'aire de $g_n(D-d_n)$ est uniform\'ement born\'ee et que la longueur de $g_n(\bo d_n)$
tend vers 0. Alors, quitte \`a extraire, $g_n$ tend au sens de Gromov vers une application holomorphe $g : D \rightarrow X$.

\fthm

En particulier
 $g_n(D_1-d_n)$ converge au sens de Hausdorff et en aire vers la r\'eunion du disque holomorphe
$g(D_1)$ et d'un nombre fini de courbes rationnelles, si $D_1$ est un disque l\'eg\`erement plus petit que $D$.

\thm{Variante 2} Soit $D^+$ le demi-disque unit\'e sup\'erieur et $g_n:D^+ \rightarrow X$ une suite
 d'applications holomorphes,
lisses jusqu'\`a $]-1,1[$. On suppose que l'aire de $g_n(D^+)$ est uniform\'ement born\'ee et
 que la longueur de
$g_n(]-1,1[)$ tend vers 0. Alors, quitte \`a extraire, $g_n$ tend au sens de Gromov vers un point.

\fthm
De m\^eme $g_n(D_1 ^+)$ converge
au sens de Hausdorff et en aire vers une union finie de courbes rationnelles si $D_1$ est un peu plus petit que $D$.

\null

Esquissons par exemple la d\'emonstration de la deuxi\`eme variante.

On peut toujours supposer que $g_n(]-1,1[)$ tend vers un point $p$. Il s'agit de v\'erifier la convergence
de $g_n$ vers $p$ au voisinage de $]-1,1[$ hors d'un nombre fini de points au bord.
Fixons $\epsilon$ petit par rapport au diam\`etre d'un domaine de carte de $X$ centr\'e en $p$.
Un point de $]-1,1[$ est dit {\it d'explosion} si, pour tout demi-disque $d^+$ centr\'e en ce point,  on a
$\overline { \lim} \text { aire}(g_n(d^+)) > \epsilon^2 $.
Soit un demi-disque $d^+$ centr\'e sur $]-1,1[$ hors de ces points d'explosion, assez petit pour que l'aire de
$g_n(d^+)$ soit born\'ee par $\epsilon^2$. Il suffit de montrer
 que $g_n$ converge vers $p$ sur $d^+$.
Par un argument longueur-aire et quitte \`a r\'eduire un peu $d^+$, on peut supposer que la longueur de
 $g_n(\bo d^+ \cap (Im(z)>0))$ est de
l'ordre de $\epsilon$. Autrement dit $g_n(\bo d^+)$ reste pr\`es de $p$. Donc tout le disque $g_n( d^+)$ reste
 dans le domaine de carte,
 sinon son aire serait
sup\'erieure \`a $\epsilon^2$ par le th\'eor\`eme de Lelong.
 Soit $h_n$ la compos\'ee de $g_n$ avec la carte sur $d^+$.
Elle est uniform\'ement born\'ee sur $d^+$ et tend vers 0 sur
 $d^+ \cap ]-1,1[$. Par le th\'eor\`eme des deux constantes $\log \Vert h_n\Vert$ tend
vers $-\infty$, d'o\`u le r\'esultat.

\null

Analysons l'apparition de bulles en un point d'explosion au bord.
On va convertir celui-ci en une explosion int\'erieure par reparam\'etrage, en raisonnant 
par induction sur la borne d'aire.
 Amor\c cons cette r\'ecurrence quand la borne d'aire vaut 1.
Consid\'erons une explosion au bord, par exemple en 0. On peut trouver une suite
$a_n$ tendant vers 0 telle que $g_n(a_n)$ converge vers un point diff\'erent de $p$. Soit $d_n^+$
le demi-disque centr\'e en $Re(a_n)$ de rayon $2Im(a_n)$ et $r_n$ son param\'etrage affine par $D^+$.
La compos\'ee $g_n\circ r_n$ satisfait la m\^eme borne d'aire et tend vers $p$ sur $]-1,1[$. Elle converge
donc vers $p$ en dehors de points d'explosion dont l'un au moins est \`a l'int\'erieur (en $\frac{i}{2}$).
 Par le th\'eor\`eme
 il s'y cr\'ee une bulle d'aire au moins 1. Elle absorbe donc toute l'aire de l'image, interdisant 
d'autres explosions pour $g_n\circ r_n$. En revenant \`a $g_n$ on obtient bien une bulle
 en 0 qui concentre l'aire.

\null

  \subhead 2. D\'emonstration du th\'eor\`eme, cas r\'eparti \endsubhead

\null

Soit $T$ un courant d'Ahlfors provenant des disques $\Delta_n$ chargeant un compact $K$. Il s'agit de construire une courbe
enti\`ere obtenue par reparam\'etrage de disques de $\Delta_n$ et passage \`a la limite au sens de Gromov et coupant $K$
sur un ensemble d'aire non nulle. Dans la suite
on parlera de courbe enti\`ere {\it issue} de $(\Delta_n)$ {\it coupant} $K$.
Notons au passage
qu'une courbe rationnelle 
produite par explosion d'une suite
de disques d'aire born\'ee dans $\Delta_n$ est un cas particulier de
courbe enti\`ere issue de $(\Delta_n)$, celui d'aire finie.

\null

Pr\'ecisons le contexte. On
 param\`etre $\Delta_n$ par une application $f_n$ holomorphe sur $D$ et lisse jusqu'au bord. Notons
$a_n$ l'aire de $\Delta_n$ et $l_n$ la longueur de son bord. Par hypoth\`ese $l_n=o(a_n)$.
Si $a_n$ reste born\'ee, $l_n$ tend vers 0. Il s'ensuit que $f_n$ tend vers une constante au sens de Gromov : toute l'aire
se concentre dans des bulles (voir pr\'eliminaires). On cr\'ee ainsi directement des courbes rationnelles issues de $(\Delta_n)$
dont l'une au moins coupera $K$.

On supposera donc que $a_n$ tend vers l'infini. Pour all\'eger l'expos\'e on fera m\^eme l'hypoth\`ese suivante :

\noindent
(H) \hskip1.4cm aucune courbe rationnelle
issue de $(\Delta_n)$ ne coupe $K$.

\null

Les courbes enti\`eres que l'on va construire seront donc d'aire infinie.

Dans la suite, une famille de disques ou d'anneaux de cardinal de l'ordre de $a_n$ sera dite
{\it consistante}. On omettra de pr\'eciser
``quitte \`a extraire en $n$'' ou ``pour $n$ assez grand'' dans les raisonnements qui vont suivre.

 \null

{\bf Construction des germes.}

\null
Comme $T$ charge $K$, les voisinages de $K$
contiennent une proportion fixe de l'aire des disques $\Delta_n$.
Autrement dit, il existe $\delta >0$ tel que
 aire$(\Delta_n \cap U_n) \geq \delta
 a_n$,
o\`u $(U_n)$ est une base de voisinages de $K$.

Notons $\mu_n$ la mesure $f_n^*({\bold 1}_{U_n}\omega)$ o\`u $\omega$ est la forme d'aire associ\'ee \`a la m\'etrique
hermitienne sur $X$. 

Construisons une famille consistante $F_n$ de disques disjoints dans $\C$ de masse 1 pour $\mu_n$.
Pour cela, \`a chaque point de $D$ on associe le plus petit disque centr\'e en ce point et de masse 1 pour $\mu_n$. 
Par le lemme de Besicovitch (voir [\MAT]), on peut extraire de cette famille
 des sous-familles en nombre au plus $C$ (constante universelle),
 chacune constitu\'ee
de disques disjoints, et dont la r\'eunion totale recouvre $D$. Il s'ensuit qu'une de ces sous-familles
 couvre un ensemble de
 masse au moins
 $\frac{\delta a_n}{C}$ pour $\mu_n$. Elle est donc consistante, c'est notre famille $F_n$. A ce stade
 les disques
de $F_n$ ne sont pas forc\'ement dans $D$.

\null

{\bf Doublement des disques.}

\null

Supposons que nous puissions extraire de $F_n$ une sous-famille $F_n^1$ consistante dont les disques doubl\'es
restent disjoints. Ici le double $2d$ (plus g\'en\'eralement $\lambda d$) d'un disque $d$ est le disque concentrique de rayon double
(de rayon multipli\'e par $\lambda$).
 On va voir qu'une bonne partie des disques
de $F_n^1$ restent dans $D$ et que leurs images sont d'aire born\'ee.

En effet, on sait que card$(F_n^1)\geq \epsilon a_n$. Un calcul d'aire montre que les disques doubl\'es de $F_n^1$ d'aire grande
(tels que aire($f_n(2d \cap D))\geq \frac{2}{\epsilon}$) sont au plus $\frac {\epsilon}{2} a_n$. Ainsi, quitte \`a r\'eduire $F_n^1$
 de moiti\'e,
les images des disques doubl\'es satisferont une borne d'aire uniforme.
De m\^eme un calcul de longueur donne que les disques doubl\'es de $F_n^1$ \`a bord long (tels que long($f_n(2d \cap \bo D))\geq
\frac{4l_n}{\epsilon a_n}$) sont au plus $\frac{\epsilon}{4}a_n$.  Donc, quitte \`a r\'eduire encore $F_n^1$,
les disques doubl\'es v\'erifieront que long($f_n(2d \cap \bo D))$ tend vers $0$ puisque $l_n=o(a_n)$.

Ceci permet de confiner les disques de $F_n^1$ dans $D$. Sinon, soit $d_n$ dans $F_n^1$ rencontrant $\bo D$. Son double
 $2d_n$ coupe largement $\bo D$. Soit $h_n$ la repr\'esentation conforme sym\'etrique de
 $D^+$ sur $2d_n \cap D$ envoyant $]-1,1[$ sur $2d_n \cap \bo D$.
Gr\^ace aux pr\'eliminaires, $f_n\circ h_n$ converge au sens de Gromov vers un point. En particulier toute
 l'aire
pr\`es de $K$ est absorb\'ee dans des bulles, contredisant (H).

\null

Continuons la discussion en s\'eparant deux cas, suivant
que l'on peut poursuivre ind\'efiniment ce processus de doublement ou non. 
Un mod\`ele simple du premier cas est une famille $F_n$ de disques de m\^eme taille bien r\'epartis dans $D$.
Le second cas correspond par exemple \`a une famille $F_n$ form\'ee d'une cha\^\i ne de disques proches
de taille d\'ecroissant rapidement : ils se concentrent en un point.
 On voit dans le cas r\'eparti qu'un nombre consistant de disques doubl\'es
 restent disjoints. Par contre, dans le cas concentr\'e les disques doubl\'es sont tous embo\^\i t\'es.

\null

{\bf Cas r\'eparti.}

\null

Finissons l'argument dans le cas o\`u l'on peut toujours doubler. On construit donc des familles consistantes
 $F_n^k$ d\'ecroissantes en $k$
 telles que :

i) si $d$ est un disque de $F_n^k$, $2^kd$ est contenu dans $D$,

ii) il existe une constante $C_k$ ind\'ependante de $n$ telle que aire($f_n(2^kd)) \leq C_k$.

On en d\'eduit ainsi une suite de disques $(d_n)$
 telle que $2^nd_n$ reste dans $D$, aire($f_n(d_n)\cap U_n)=1$ et aire($f_n(2^kd_n)) \leq C_k$ pour $1\leq k\leq n$.
Si $r_n$ est le param\'etrage affine de $d_n$ par $D$, $f_n\circ r_n$ est donc d\'efinie sur $2^nD$ et d'aire
 uniform\'ement born\'ee
sur les compacts de $\C$. Cette suite converge au sens de Gromov vers une courbe enti\`ere qui coupe $K$.

\null

  \subhead 3. D\'emonstration du th\'eor\`eme,
 cas concentr\'e
 \endsubhead

\null

Supposons que le processus de doublement s'interrompe \`a un moment, par exemple au d\'ebut.
Autrement dit, on ne peut extraire de $F_n$ une famille consistante dont les disques
doubl\'es soient disjoints.
On va profiter de ceci pour construire des lign\'ees de disques
grossis embo\^ \i t\'es, comme dans le mod\`ele concentr\'e.
Soit  $4F_n$ la famille des disques $4d$ o\`u $d$ parcourt $F_n$.

\null

{\bf Pr\'eparation de la famille. }

\null

Quitte \`a r\'eduire $4F_n$ d'un facteur fixe, on peut supposer que seuls les disques de taille
franchement diff\'erente s'intersectent. Pr\'ecis\'ement, soient deux disques $d,d'$ dans $4F_n$ ($d'$ plus petit que $d$).
 On veut s'assurer que si $2d$ et $ 2d'$ s'intersectent,
alors $d'$ est 16 fois plus petit que $d$.

Pour cela, d\'enombrons les disques $d'$ plus grands que $\frac {1}{16}d$ tels que $2d$ et $2d'$ s'intersectent.
N\'ecessairement $5d$ contient $d'$. Or, les disques $\frac{1}{4}d'$ \'etant disjoints et plus grands que $\frac {1}{64}d$,
ils sont en nombre au plus $320^2$ dans $5d$ par un calcul d'aire. Il suffit donc de r\'eduire notre famille
de ce facteur. Notons-la encore $4F_n$.

\null

{\bf Arborescence d'incidence. }

\null

Les incidences de disques donnent une structure de graphe orient\'e sur $4F_n$. Les sommets en
sont les disques et une ar\^ete joint $d$ \`a $d'$ si $d'$ est plus petit que $d$ et $d'$ intersecte
$d$.

Voici comment on modifie ce graphe en une arborescence.

 Partons du disque $d$ le plus petit de $4F_n$. S'il re\c coit $k$ ar\^etes, leurs origines $d_1,...,d_k$
 sont de taille tr\`es diff\'erente par ce qui pr\'ec\`ede. On les suppose ordonn\'ees de mani\`ere d\'ecroissante.
On remplace alors les ar\^etes d'extr\'emit\'e $d$ par le chemin d'ar\^etes $d_1 \rightarrow d_2 \rightarrow
...\rightarrow d_k \rightarrow
d$. On proc\`ede de m\^eme avec le disque suivant en taille en tenant compte des
modifications pr\'ec\'edentes et l'on continue jusqu'\`a \'epuisement de la famille.

On obtient une arborescence qui code encore les (presque) incidences des disques. En effet, si $d$ intersecte
$d'$ on a toujours un chemin dans l'arborescence de $d$ \`a $d'$. Inversement, en pr\'esence d'un
chemin de $d$ \`a $d'$ on v\'erifie que $d'$ est contenu dans $2d$.

\null

{\bf Analyse de l'arborescence. }

\null

Introduisons un peu de terminologie. 
Appelons {\it fils} de $d$ l'extr\'emit\'e d'une ar\^ete issue de $d$ (on parlera inversement
de {\it p\`ere}).

On distingue dans l'arborescence les sommets {\it terminaux} (sans fils), {\it simples}
(\`a fils unique) et {\it multiples} (avec au moins deux fils). Notons
$t_n$, $s_n$ et $m_n$ leurs nombres respectifs. Nous allons voir que la plupart des sommets sont simples.

 Le nombre total de sommets est de l'ordre de $a_n$.
Par construction les sommets terminaux correspondent \`a des disques
disjoints entre eux dans $4F_n$. Ils sont donc par hypoth\`ese en nombre n\'egligeable devant $a_n$.
En d\'enombrant de deux mani\`eres les ar\^etes (par leur origine ou leur extr\'emit\'e), on
v\'erifie classiquement que la somme des valences est major\'ee par le nombre de sommets.

Autrement dit,
$ s_n +2m_n \leq t_n+s_n+m_n $
d'o\`u $m_n \leq t_n$. Donc $m_n+t_n = o(a_n)$, et la plupart des sommets sont simples.

Supprimons les sommets multiples de l'arborescence. Celle-ci se d\'ecompose alors en un certain
nombre de chemins sans branchement appel\'es {\it lign\'ees}.
Leur nombre est n\'egligeable devant $a_n$ puisque toute lign\'ee aboutit \`a un sommet terminal ou
multiple.

\null

{\bf Pr\'eparation des lign\'ees. }

\null

Ces lign\'ees correspondent essentiellement \`a des successions de disques embo\^\i t\'es.
Pr\'ecis\'ement, $2d'$ est contenu dans $\frac {3}{4}d$ pour la plupart des ar\^etes $d \rightarrow d'$ des lign\'ees.

Si ce n'est pas le cas, $2d'$ \'evite $\frac{1}{2}d$ puisque $d'$ est beaucoup plus petit
que $d$. Donc
tous les descendants de $d'$ sont disjoints 
de $\frac{1}{2}d$ par ce qui pr\'ec\`ede. Autrement dit, le disque $ \frac{1}{2}d$ est un des sommets
 terminaux de la famille $2F_n$ munie
 de son graphe
d'incidence. Comme ceux-ci sont disjoints, leur nombre est n\'egligeable devant $a_n$.

 Supprimons ces exceptions. Cela coupe nos lign\'ees en des lign\'ees plus courtes, en nombre
toujours n\'egligeable devant $a_n$. Les anneaux obtenus comme diff\'erence de deux disques d'une m\^eme lign\'ee
 vont jouer dans le cas concentr\'e le r\^ole d\'evolu
aux disques dans le cas r\'eparti.

Remarquons d\'ej\`a que l'image d'un anneau $d_1-d_7$, o\`u $d_1\rightarrow ... \rightarrow d_7$ est
 un morceau de lign\'ee
de longueur 6, conserve de l'aire pr\`es de $K$. En effet, cet anneau (et m\^eme l'anneau plus petit
 $\frac{3}{4}d_1-2d_7$) doit
contenir l'un des disques initiaux $\frac {1}{4}d_1,...,\frac {1}{4}d_6$ de
 $F_n$. Sinon ces disques
intersecteraient tous $2d_7$. On pourrait alors les grossir l\'eg\`erement en des disques d'intersection non vide.
C'est g\'eom\'etriquement impossible puisque chaque disque est centr\'e en dehors des autres (voir [\MAT], p.29).

\null

Pour all\' eger les notations, on \'eclaircit nos lign\'ees de la mani\`ere suivante : on les subdivise en morceaux
de longueur 6, puis on remplace chaque morceau par une seule ar\^ete en supprimant les sommets interm\'ediaires.

Pour ces nouvelles lign\'ees, les anneaux de la forme $d-d'$ ($d'$ fils de $d$)
sont par construction disjoints, de module minor\'e et d'image conservant de l'aire
pr\`es de $K$. Ils forment une famille consistante not\'ee $G_n$.

\null

{\bf Doublement des anneaux. }

\null  

Voici comment on ``double'' ces anneaux. Si $a$ est l'anneau $d-d'$, on notera $2a$ l'anneau $d_1-d_1'$ o\`u
$d_1$ est le p\`ere de $d$ et $d_1'$ le fils de $d'$.
On peut doubler la plupart des anneaux de $G_n$ et ces doubles se recouvrent au plus 3 fois.
On peut donc, comme au paragraphe 2, extraire de $G_n$
une famille consistante $G^1_n$ dont les anneaux doubl\'es sont disjoints et v\'erifient :

i)  aire($f_n(2a\cap D)) \leq C_1$,

ii) long($f_n(2a\cap \partial D))$ tend vers 0.

\null

Comme pour les disques dans le cas r\'eparti, ceci force les anneaux de 
$G_n^1$ \`a rester dans $D$.
 Sinon, soit $a_n=d-d'$ dans $G_n^1$ rencontrant $\bo D$. Son double $d_1-d_1'$ coupe largement $\bo
D$. En fait $\bo D$ s\'epare l'anneau ext\'erieur $b_n=d_1-d$ de son double
en deux disques topologiques. Soit $h_n$ la repr\'esentation conforme sym\'etrique de $D^+$ sur $d_1\cap D$ envoyant $]-1,1[$ sur
$d_1\cap \bo  D$.
La suite diff\`ere suivant que $d$ et $d_1$ restent de taille
comparable ou non.

Dans le premier cas, la pr\'eimage par $h_n$ de $b_n \cap D$
converge vers un disque topologique dont le bord
rencontre largement $]-1,1[$. Gr\^ace aux pr\'eliminaires,
 $f_n \circ h_n$ y tend au sens de Gromov vers un point. Toute l'aire de l'image pr\`es de $K$
 doit \^etre absorb\'ee dans des bulles, contredisant (H).

Dans le second cas, le module de $b_n$ tend vers l'infini. Par un argument longueur-aire, on trouve
un cercle $c_n$ dans $b_n$ tel que long($f_n(c_n\cap D))$ tend vers 0. Ce cercle s\'epare $b_n$ en
un anneau ext\'erieur $b_n^+$ et un anneau int\'erieur $b_n^-$ dont les modules tendent aussi
vers l'infini. L'image de l'un d'entre eux, par exemple $b_n^+$,
conserve de l'aire pr\`es de $K$. La
pr\'eimage de $b_n^+ \cap D$ par $h_n$ tend vers $D^+$. Comme plus haut $f_n \circ h_n$ y converge au sens de Gromov vers
un point. Ceci force l'apparition de courbes rationnelles coupant $K$, contredisant (H).

\null

{\bf Fin de l'argument. }

\null

Cette fois, le processus de doublement s'it\`ere sans obstacle.
 On construit donc des familles
consistantes $G_n^k$ d\'ecroissantes en $k$ telles que :

i) si $a$ est un anneau de $G_n^k$, $2^ka$ est bien d\'efini et contenu dans $D$,

ii) il existe une constante $C_k$ ind\'ependante de $n$ telle que aire($f_n(2^ka))\leq C_k$.

On en d\'eduit ainsi une suite d'anneaux $(a_n)$
 telle que $2^na_n$ reste dans $D$, aire($f_n(a_n)\cap U_n)\geq 1$
et aire($f_n(2^ka_n)) \leq C_{k}$ pour $1\leq k\leq n$.

 Comme plus haut, la fin de l'argument diff\`ere suivant que les bords int\'erieur et ext\'erieur de $a_n$ restent
de taille comparable ou non.

Dans le premier cas, soit
$r_n$ l'automorphisme de $\C$ envoyant l'origine
sur le centre du bord int\'erieur de $2^na_n$ et de rapport le diam\`etre de $a_n$.
 Par construction la pr\'eimage de $2^ka_n$ par $r_n$ tend vers $\C^*$ 
quand $k$ tend vers l'infini. La compos\'ee $f_n\circ r_n$ est donc d'aire uniform\'ement born\'ee sur les compacts de $\C^*$.
Elle converge au sens de Gromov vers $f:\C^* \rightarrow X$.
Comme l'anneau central $r_n^{-1}(a_n)$ reste dans un compact fixe de $\C^*$, on en conclut 
 que $f(\C^*)$ coupe $K$ puisque l'aire pr\`es de $K$ ne peut \^etre absorb\'ee dans des bulles d'explosion.

Dans le second cas, le module de $a_n$ tend vers l'infini et l'on peut trouver un cercle
dans $a_n$ dont l'image par $f_n$ a une longueur qui tend vers 0.
Ce cercle s\'epare $a_n$
 en deux anneaux $a^+_n$ et $a^-_n$
dont le module tend aussi vers l'infini.
L'image de l'un d'entre eux par $f_n$, par exemple $a^+_n$, conserve de l'aire
pr\`es de $K$.
Soit $r_n$ l'automorphisme de $\C$ envoyant l'origine sur le centre du bord int\'erieur de $a_n^+$
et de rapport le diam\`etre de $a_n$. La pr\'eimage de $2^k
a^+_n$ par $r_n$ tend vers $\C^*$ quand $k$ tend vers l'infini. Gr\^ace aux pr\'eliminaires,
$f_n\circ r_n$ converge au sens de Gromov vers une courbe enti\`ere qui coupe $K$.

 \Refs

\widestnumber\no{99}
\refno=0

\bref \by R. Brody \paper Compact manifolds and hyperbolicity \jour Trans. Amer. Math. Soc. \vol235\yr1978\pages213--219
\endref

\bref \by M. Brunella \paper Courbes enti\`eres et feuilletages holomorphes \jour
Ens. Math. \vol45\yr1999\pages195--216
\endref

\bref \by J. Duval \paper Singularit\'es des courants d'Ahlfors \jour Ann. Sci. ENS \vol39\yr2006\pages527--533
\endref

\bref \by  B. Kleiner \paper Hyperbolicity using minimal surfaces \jour preprint
\endref

\bref \by M. McQuillan \paper Diophantine approximations and foliations \jour Publ. Math. IHES \vol87\yr1998\pages121--174
\endref

\bref \by M. McQuillan \paper Bloch hyperbolicity \jour preprint IHES 2001
\endref

\bref \by P. Mattila \book Geometry of sets and measures in euclidean spaces\publ Cambridge University Press \yr 1995 \publaddr Cambridge
\endref

\bref \by M. Paun \paper Currents associated to transcendental entire curves on compact K\"ahler manifolds \jour preprint 2003
\endref

\bref \by J.-C. Sikorav \book Some properties of holomorphic curves in almost complex manifolds, {\rm in} Holomorphic curves in symplectic geometry
\pages165--189 \bookinfo Prog. Math. \vol117 \publ Birkh\"auser \yr1994 \publaddr Basel
\endref
\endRefs

\address 
\noindent  
Laboratoire \'Emile Picard, 
  Universit\'e Paul Sabatier, 31062 Toulouse Cedex 09.
 \endaddress
\email 
  duval\@picard.ups-tlse.fr 
\endemail

\enddocument